\input amstex
\documentstyle{amsppt}
\magnification=\magstephalf
 \addto\tenpoint{\baselineskip 15pt
  \abovedisplayskip18pt plus4.5pt minus9pt
  \belowdisplayskip\abovedisplayskip
  \abovedisplayshortskip0pt plus4.5pt
  \belowdisplayshortskip10.5pt plus4.5pt minus6pt}\tenpoint
\pagewidth{6.5truein} \pageheight{8.9truein}
\subheadskip\bigskipamount
\belowheadskip\bigskipamount
\aboveheadskip=3\bigskipamount
\catcode`\@=11
\def\output@{\shipout\vbox{%
 \ifrunheads@ \makeheadline \pagebody
       \else \pagebody \fi \makefootline 
 }%
 \advancepageno \ifnum\outputpenalty>-\@MM\else\dosupereject\fi}
\outer\def\subhead#1\endsubhead{\par\penaltyandskip@{-100}\subheadskip
  \noindent{\subheadfont@\ignorespaces#1\unskip\endgraf}\removelastskip
  \nobreak\medskip\noindent}
\outer\def\enddocument{\par
  \add@missing\endRefs
  \add@missing\endroster \add@missing\endproclaim
  \add@missing\enddefinition
  \add@missing\enddemo \add@missing\endremark \add@missing\endexample
 \ifmonograph@ 
 \else
 \vfill
 \nobreak
 \thetranslator@
 \count@\z@ \loop\ifnum\count@<\addresscount@\advance\count@\@ne
 \csname address\number\count@\endcsname
 \csname email\number\count@\endcsname
 \repeat
\fi
 \supereject\end}
\catcode`\@=\active
\CenteredTagsOnSplits
\NoBlackBoxes
\nologo
\def\today{\ifcase\month\or
 January\or February\or March\or April\or May\or June\or
 July\or August\or September\or October\or November\or December\fi
 \space\number\day, \number\year}
\define\({\left(}
\define\){\right)}

\define\CC{{\Bbb C}}
\define\CP{{\Bbb C\Bbb P}}

\define\End{\operatorname{End}}

\define\Hom{\operatorname{Hom}}

\define\RR{{\Bbb R}}

\define\Tr{\operatorname{Tr}}
\define\ZZ{{\Bbb Z}}
\define\[{\left[}
\define\]{\right]}

\define\chiup{\raise.5ex\hbox{$\chi$}}
\define\cir{S^1}

\define\dbar{{\bar\partial}}

\define\exertag #1#2{#2\ #1}

\define\inv{^{-1}}
\define\mstrut{^{\vphantom{1*\prime y}}}
\define\protag#1 #2{#2\ #1}
\define\rank{\operatorname{rank}}
\define\res#1{\negmedspace\bigm|_{#1}}
\define\temsquare{\raise3.5pt\hbox{\boxed{ }}}

\define\theprotag#1 #2{#2~#1}

\define\xca#1{\removelastskip\medskip\noindent{\smc%
#1\unskip.}\enspace\ignorespaces }

\define\zmod#1{\ZZ/#1\ZZ}

\define\zt{\zmod2}

\define\Atil{\tilde{A}}

\define\Fred{\operatorname{Fred}}
\define\GR{G(R)}
\define\Gtil{\tilde{G}}
\define\Gt{\tilde{\G}}
\define\G{\Cal{G}}
\define\Pin{\operatorname{Pin}}
\define\TT{\Bbb{T}}
\define\fatil{\tilde{\frak{a}}}
\define\fa{\frak{a}}
\define\fc{\scrF_C}
\define\hol{\operatorname{hol}}
\define\id{\operatorname{id}}
\define\lgc{LG^{\sigma _R(G)}_R}
\define\lgt{LG^{\tau }_R}
\define\ring{\RR[[u,u\inv ]]}
\define\scrF{\Cal{F}}

\define\ttil{\tilde\tau }
\define\zc{Z_C}
\redefine\lg{LG\mstrut _R}

\NoRunningHeads

\refstyle{A}
\widestnumber\key{SSSSS}   
\document

	\topmatter
 \title\nofrills{Twisted $K$-Theory and Loop Groups}\endtitle 
 \thanks The author is supported by NSF grant DMS-0072675.\endthanks
 \affil Department of Mathematics \\ University of Texas at Austin\endaffil 
 \address Department of Mathematics, University of Texas, Austin, TX
78712\endaddress 
 \email dafr\@math.utexas.edu \endemail
 \author Daniel S. Freed  \endauthor
 \date June 23, 2002\enddate
 \abstract
 Twisted $K$-theory has received much attention recently in both mathematics
and physics.  We describe some models of twisted $K$-theory, both topological
and geometric.  Then we state a theorem which relates representations of loop
groups to twisted equivariant $K$-theory.  This is joint work with Michael
Hopkins and Constantin Teleman.
	\endabstract
	\endtopmatter

\document

 \comment
 lasteqno @ 12
 \endcomment

The loop group of a compact Lie group~$G$ is the space of smooth maps $S^1\to
G$ with multiplication defined pointwise.  Loop groups have been around in
topology for quite some time~\cite{Bo}, and in the 1980s were extensively
studied from the point of view of representation theory~\cite{Ka}, \cite{PS}.
In part this was driven by the relationship to conformal field theory.  The
interesting representations of loop groups are projective, and with fixed
projective cocycle~$\tau $ there is a finite number of irreducible
representations up to isomorphism.  Considerations from conformal field
theory~\cite{V} led to a ring structure on the abelian group~$R^\tau (G)$
they generate, at least for {\it transgressed\/} twistings.  This is the {\it
Verlinde ring\/}.  For $G$~simply connected $R^\tau (G)$~is a quotient of
the representation ring of~$G$, but that is not true in general.  At about
this time Witten~\cite{W} introduced a three-dimensional topological quantum
field theory in which the Verlinde ring plays an important role.  Eventually
it was understood that the fundamental object in that theory is a ``modular
tensor category'' whose Grothendieck group is the Verlinde ring.  Typically
it is a category of representations of a loop group or quantum group.
 
For the special case of a finite group~$G$ the topological field theory is
specified by a certain cocycle on~$G$ and the category can be calculated
explicitly~\cite{F1}.  We identified it as a category of representations of a
Hopf algebra constructed from~$G$, thus directly linking the Chern-Simons
lagrangian and quantum groups.  Only recently did we realize that this
category has a description in terms of {\it twisted\/} equivariant
$K$-theory, and it was natural to guess that the Verlinde ring for
arbitrary~$G$ has a similar description.  Ongoing joint work with Michael
Hopkins and Constantin Teleman has confirmed this description.  We can
speculate further and hope that twisted $K$-theory provides a construction of
the modular tensor category, and perhaps even more of the three-dimensional
topological field theory.  In another direction the use of $K$-theory may
shed light on Verlinde's formula for certain Riemann-Roch numbers.  In any
case our result fits well with other uses of $K$-theory in representation
theory~\cite{CG}, for example in the geometric Langlands program.
 
The physics motivation for the main theorem is discussed in~\cite{F2},
\cite{F3}.  Here, in~\S{3}, we explain the statement of our result in
mathematical terms; the proof will appear elsewhere.  As background we
describe some concrete topological models of twisted $K$-theory in~\S{1}, and
give a twisted version of the Chern-Weil construction in~\S{2}.

As mentioned above the work I am discussing is being carried out with Michael
Hopkins and Constantin Teleman.  I thank them for a most enjoyable
collaboration.

 \head
 \S{1} Twistings of $K$-theory
 \endhead

Let $X$~be a reasonable compact space, say a finite CW~complex.  Then
isomorphism classes of complex vector bundles over~$X$ form a semigroup whose
group completion is the $K$-theory group~$K^0(X)$.  This basic idea was
introduced by Grothendieck in the context of algebraic geometry~\cite{BS},
and was subsequently transported to topology by Atiyah and
Hirzebruch~\cite{AH}.  Vector bundles are local---they can be cut and
glued---and in the topological realm this leads to a cohomology theory.  In
particular, there are groups~$K^n(X)$ defined for~$n\in \ZZ$.  Historically,
$K$-theory was the first example of a generalized cohomology theory, and it
retains the features of ordinary cohomology with one notable exception: the
cohomology of a point is nontrivial in all even degrees, as determined by
Bott periodicity.  There are many nice spaces~$B$ which represent complex
$K$-theory in the sense that $K^0(X)$~is the set of homotopy classes of maps
from~$X$ to~$B$.  A particularly nice choice~\cite{A1,Appendix}, \cite{J} is
the space~$B=\Fred(H)$ of Fredholm operators on a separable complex Hilbert
space~$H$.  Thus a map $T\:X\to \Fred(H)$ determines a $K$-theory class
on~$X$.  It is convenient to generalize and allow the Hilbert space~$H$ to
vary as follows.  A {\it Fredholm complex\/}\footnote{ One~\cite{S1} can
allow more general topological vector spaces and complexes which are nonzero
in degrees other than~0 and~1.} is a graded Hilbert space bundle
$E^\bullet=E^0 \oplus E^1\to X$ together with a fiberwise Fredholm
map~$E^0\overset T \to \to E^1$, and it also represents an element
of~$K^0(X)$.  Another innovation was the introduction by Atiyah and
Segal~\cite{S2} of the equivariant $K$-theory groups~$K^\bullet_G(X)$ for a
compact space~$X$ which carries the action of a compact Lie group~$G$.  The
basic objects are $G$-equivariant vector bundles~$E\to X$, and $K^0_G(X)$~is
the group completion of the set of equivalence classes.  For example, if
$X$~is a point then the equivariant $K$-theory is the representation
ring~$K_G$ of the compact Lie group~$G$; in general, $K_G(X)$~is a
$K_G$-module.
 
As a first example of twisted $K$-theory we consider twisted versions
of~$K_G$.  A twisting~$\tau $ is a central extension 
  $$ 1 \longrightarrow \TT \longrightarrow \Gtil \longrightarrow
     G\longrightarrow 1,  $$
where $\TT\subset \CC$~is the circle group of unit norm complex numbers.
Then the twisted representation ``ring''~$K^\tau _G$ is the group completion
of equivalence classes of complex representations of~$\Gtil$ on which the
central~$\TT$ acts by scalar multiplication.  Twisted $K$-theory is {\it
not\/} a ring, but rather $K^\tau _G$ is a $K_G$-module.  A map of twistings,
i.e., an isomorphism of central extensions, determines an isomorphism of
twisted $K$-groups.  So twisted $K$-theory is determined up to {\it
noncanonical\/} isomorphism by the equivalence class of the twisting.  For
example, $K_{SO(3)}\cong \ZZ[s]$ is the polynomial ring on a single
generator, the 3-dimensional defining representation.  Up to equivalence
there is a single nontrivial central extension
  $$ \tau =\{\;1 \longrightarrow \TT \longrightarrow U(2) \longrightarrow
     SO(3)\longrightarrow 1\;\}  $$
which is induced by the inclusion $\zt\hookrightarrow \TT$ from the extension
  $$ 1 \longrightarrow \zt \longrightarrow SU(2) \longrightarrow
     SO(3)\longrightarrow 1.  $$
Virtual representations of~$U(2)$ on which the center acts naturally
correspond 1:1 with virtual representations of~$SU(2)$ on which the central
element acts as~$-1$.  Now~$K_{SU(2)}\cong \ZZ[t]$, where $t$~is the defining
2-dimensional representation, and we identify~$K^\tau _{SO(3)}\subset
K_{SU(2)}$ as the subgroup of odd polynomials in~$t$; the $K_{SO(3)}\cong
\ZZ[s]$-module structure is~$s\cdot t=t^3-t$, by the Clebsch-Gordon rule, and
$K^\tau _{SO(3)}$~is a free module of rank one.
 
More generally, on a $G$-space~$X$ a ``cocycle''~$\tau $ for the equivariant
cohomology group~$H^3_G(X;\ZZ)$ defines a twisted $K$-theory group~$K^\tau
_G(X)$ which is a module over~$K^0_G(X)$.  A better point of view is that
$\tau $~is a cocycle, or geometric representative, of a class
in~$H^1_G(X;A_G)$, where $A_G$~is a group of automorphisms of equivariant
$K$-theory.  We will not try to make ``automorphism of equivariant
$K$-theory'' precise here, but content ourselves with a concrete model, first
in the nonequivariant case.  Take~$\Fred(H)$ to be the classifying space of
$K$-theory.  Then the group~$A=PGL(H)$ acts as automorphisms by conjugation.
Since $\Atil=GL(H)$~is contractible~\cite{K}, the
quotient~$A=\Atil/\CC^{\times }$ is homotopy equivalent to~$\CP^{\infty}$,
and so~$H^1(X;A)\cong H^3(X;\ZZ)$.  A twisting~$\tau $ can be taken to be a
principal~$A$-bundle $\pi \:P\to X$.  The action of~$A$ on~$\Fred(H)$ defines
an associated bundle~$\Fred(H)_P\to X$, and the twisted $K$-group~$K^\tau
(X)$~ is the group of homotopy classes of sections of this bundle~\cite{A2}.
A section is an $A$-equivariant map $T\:P\to \Fred(H)$, and as before it is
convenient to let the Hilbert space vary.  Thus define a {\it $P$-twisted
Fredholm complex\/} to be an $\Atil$-equivariant graded Hilbert space
bundle~$E^\bullet\to P$, where the center~$\CC^\times\subset \Atil$ acts by
scalar multiplication, together with a fiberwise $\Atil$-equivariant Fredholm
operator $E^0\overset T \to \to E^1$.  Then $(E^\bullet,T)$~represents an
element of~$K^\tau (X)$.

It is perhaps unsettling that the model is infinite dimensional, but that is
unavoidable unless the class of the twisting in~$H^3(X;\ZZ)$ is torsion.
 
There are many other models of twistings and twisted $K$-theory.  For
example, {\it gerbes\/} are geometric representatives of elements in degree
three integral cohomology; see~\cite{H}, \cite{B}, \cite{M} for example.  In
a \v Cech description we have a covering~$X=\bigcup_i U_i$ of~$X$ by open
sets and a complex line bundle~$L_{ij}\to U_i\cap U_j$ on double
intersections.  There is further cocycle data on triple intersections.  In a
model of~$K$-theory on which line bundles act as automorphisms this can be
used to define twisted $K$-theory; see~\cite{BCMMS} for example.  In fact,
the ``group''~$\zt\times \CP^{\infty}$ of {\it graded\/} lines act as
automorphisms of $K$-theory, so there is a larger group of (equivalence
classes of) twistings
  $$ H^1(X;\zt \times \CP^{\infty})\cong H^1(X;\zt)\times H^3(X;\ZZ).
     \tag{1} $$
We remark that there is a natural group structure on~\thetag{1}, but it is
not the product.  In topology these twisted versions of $K$-theory, at least
for torsion twistings, were introduced by Donovan and Karoubi~\cite{DK}, who
also considered real versions.  There is another viewpoint and generalization
of $K$-theory using $C^*$-algebras, and in that context twisted $K$-theory
was discussed by Rosenberg~\cite{R} for both torsion and nontorsion
twistings.  Twisted versions of $K$-theory have appeared recently in various
parts of geometry and index theory, for example in~\cite{LU}, \cite{AR},
\cite{To}, \cite{NT}, \cite{MMS}.
 
A generalization of the previous model is useful.  Here a twisting is a
quartet~$\tau =(\G,\epsilon ,\Gt, P)$:
  $$ \alignedat2
      &\G&&\qquad \text{topological group} \\
      &\epsilon \:\G\to\zt&&\qquad \text{homomorphism termed the {\it
     grading\/}} \\
      &\Gt\to\G&&\qquad \text{central extension by $\TT$} \\
      &P\to X&&\qquad \text{principal $\G$-bundle}\endaligned \tag{2} $$
We require the existence of a homomorphism $\Gt\to GL(H)$ which is the
identity on the central~$\TT$.  Let $\Gt_0\to\G_0$ be the restriction of the
central extension over the subgroup~$\G_0=\epsilon \inv (0)$.  Then the
equivalence class of~$\tau $ is the obstruction in~\thetag{1} to
restricting/lifting~$P$ to a principal~$\G_0$ bundle.  Our previous
construction is the special case~$G=PGL(H)$ and $\epsilon $~trivial.  An
element of~$K^\tau (X)$ is represented by a $\Gt$-equivariant $\zt$-graded
Fredholm complex over~$P$, where the action of~$\Gt$ is compatible with the
grading~$\epsilon $.  If a compact Lie group~$G$ acts on~$X$ then it is easy
to extend this to a model of equivariant twistings and equivariant
$K$-theory.

As an illustration of a nontrivial $H^1$-twisting, consider~$X=pt$ and~$\tau
=\bigl(O(2),\epsilon ,O(2)\times \TT,O(2) \bigr)$ with nontrivial
grading~$\epsilon $.  The representation ring of~$O(2)$ may be written
  $$ K_{O(2)} \cong \ZZ[\sigma ,\delta ]\bigm/ (\sigma (\delta -1), \delta
     ^2-1),  $$
where $\delta $~is the one-dimensional sign representation and $\sigma $~the
standard two dimensional representation.  Then the twisted $K$-group $K^\tau
_{O(2)}$~is a module over~$K_{O(2)}$ with a single generator~$t$ and the
relation~$\delta \cdot t=t$.  There is a nontrivial odd twisted
$K$-group~$K^{\tau +1}_{O(2)}$ which is also a module with a single
generator~$u$; the relations are $\delta \cdot u=-u$ and~$\sigma \cdot u=0$.

Many topological properties of $K$-theory, including exact and spectral
sequences, have straightforward analogs in the twisted case.  This is easiest
to see from the homotopy-theoretic view of cohomology theories, and so
applies to twisted cohomology theories in general.  Computations are usually
based on these sequences.  A more specialized result is the completion
theorem in equivariant $K$-theory~\cite{AS}; its generalization to the
twisted case has some new features~\cite{FHT}.  The Thom isomorphism theorem
fits naturally into the twisted theory~\cite{DK}.  Let $V\to X$ be a real
vector bundle of finite rank, which for convenience we suppose endowed with a
metric.  There is an associated twisting
  $$ \tau (V)=\bigl(O(n),\epsilon ,\Pin^c(n),O(V)\bigr),  $$
where $O(V)\to X$ is the orthonormal frame bundle of~$V$ and $\epsilon
\:O(n)\to \zt$ is the nontrivial grading.  The isomorphism class of~$\tau
(V)$ is the pair of Stiefel-Whitney classes~$\bigl(w_1(V),W_3(V) \bigr)$.
Denote
  $$ K^{V}(X) = K^{\rank(V)+\tau (V)}(X) \tag{3} $$
as the degree-shifted twisted $K$-theory.  Then the Thom isomorphism theorem
asserts that the natural map
  $$ K^{q+V}(X) \longrightarrow K^{q+2\rank(V)}(V) \tag{4} $$
is an isomorphism.  We remark that if $X$~is a smooth compact manifold, then
Poincar\'e duality identifies $K^{TX}(X)$ with the $K$-homology
group~$K_0(X)$.  Also, we can define~$K^V(X)$ for virtual bundles~$V$.   
 
The Chern character maps twisted $K$-theory to a twisted version of real
cohomology, as we explain next in the context of Chern-Weil theory.

\head \S2 A differential-geometric model
\endhead

For simplicity we work in the nonequivariant context.  Our thoughts here were
stimulated by reading~\cite{BCMMS}.

Let $\fa=\frak{p}\frak{g}\frak{l}(H)$ denote the Lie algebra of~$\Atil=PGL(H)$
and~$\fatil=\frak{g}\frak{l}(H)$ the Lie algebra of~$A=GL(H)$.  They fit into
the exact sequence of Lie algebras
  $$ 0 \longrightarrow \CC \longrightarrow \fatil \longrightarrow \fa
     \longrightarrow 0.  $$
A linear map~$L\:\fa\to\fatil$ is a splitting if the composition with
$\fatil\to\fa$ is the identity.  A splitting~$L$ determines a closed
right-invariant 2-form on~$A$ which represents the generator of~$H^2(A;\ZZ)$;
its value on right-invariant vector fields~$\xi ,\eta $ is
  $$ \bigl[ L (\xi ),L (\eta )\bigr] - L \bigl([\xi ,\eta ]\bigr).  $$
 
Let $\pi \:P\to X$ be a principal $A$-bundle.  We now add two pieces of
geometric data: 
  $$ \alignedat{2}
      &\Theta \in \Omega ^1(P;\fa)\qquad &&\text{connection on $P\to X$} \\
      &L\:P\to\Hom(\fa,\fatil)\qquad &&\text{$A$-invariant map into
     splittings}\endaligned \tag{5} $$
Both are sections of affine space bundles over~$X$, so can be constructed
using partitions of unity.  As usual, define the curvature
  $$ F_\Theta  = d\Theta  + \frac 12[\Theta \wedge \Theta ].  $$
It is an $\fa$-valued 2-form on~$P$.  Introduce the {\it scalar\/} 2-form
  $$ \beta = \bigl(dL(\Theta )+\frac 12[L(\Theta )\wedge L(\Theta )] \bigr) -
     L(F_A).  $$
Then one can check that $\beta $~is transgressive.  In other words, $d\beta
=\pi ^*\eta $ for a closed scalar 3-form~$\eta \in \Omega ^3(X)$.  The de Rham
cohomology class of~$\eta $ in~$H^3(X;\RR)$ represents the image in real
cohomology of the isomorphism class of the twisting~$P\to X$.   
 
Now let $E^0 \overset T \to \to E^1$ be a twisted Fredholm complex over~$P$.
Thus $E^i\to P$ are $\Atil$-equivariant Hilbert space bundles, with
$\CC^\times \subset \Atil$ acting by scalar multiplication, and $T$~is an
$\Atil$-equivariant Fredholm map.  The $\Atil$~action determines an
$\Atil$-invariant partial covariant derivative on~$E^\bullet\to P$ along the
fibers of~$\pi \:P\to X$.  Again we introduce differential-geometric data:
  $$ \nabla ^\bullet\qquad \text{$\Atil$-invariant extension of the partial
     covariant derivative}  $$
Such an extension is a section of an affine space bundle over~$X$, so can be
constructed via a partition of unity.  Introduce a formal parameter~$u$ of
degree~2 and its inverse~$u\inv $ of degree~$-2$, and so the graded ring of
$\RR[[u,u\inv ]]$-valued differential forms.  (One can identify~$u$ as a
generator of the $K$-theory of a point.)  Following Quillen~\cite{Q} define
the $\Atil$-invariant superconnection
  $$ D = \pmatrix \nabla ^0&&uT^*\\T&&\nabla ^1 \endpmatrix  $$
on~$E^\bullet\to X$.  Its usual curvature $D^2\in \Omega \bigl(P;\End
E\otimes \ring\bigr)^2$ is an $A$-invariant form of total degree~2.  However,
it does not descend to the base~$X$.  Instead, one can check that the {\it
twisted\/} curvature
  $$ F(E^\bullet,T,\nabla ^\bullet) := D^2 - \beta \cdot \id  $$
is $A$-invariant and basic, so descends to an element of~$\Omega \bigl(X;\End
E\otimes \ring\bigr)^2$.  (Note that $\End E^\bullet\to P$ descends to a
graded vector bundle on~$X$, since the center of~$\Atil$ acts trivially.)  It
does not, however, satisfy a Bianchi identity, since $\beta $~is not closed.
 
The Chern character form
  $$ ch(E^\bullet,T,\nabla ^\bullet) := \Tr\exp(u\inv F) \tag{6} $$
is an element in~$\Omega \bigl(X; \ring\bigr)^0$ of total degree~0.  Here we
assume favorable circumstances in which the graded trace~$\Tr$ is finite.
For example, if the twisting class is torsion then we can take $E^\bullet\to
P$ finite dimensional.  Or, if the superconnection comes from a family of
elliptic operators as in~\cite{Bi} then the graded trace exists.  The Chern
character form~\thetag{6} is not closed in the usual sense, but rather
  $$ (d + u\inv \eta )\,ch(E^\bullet,T,\nabla ^\bullet)=0.  $$
The differential $d+u\inv \eta $ on~$\Omega \bigl(X; \ring\bigr)^\bullet$
computes a twisted version of real cohomology which is the codomain of the
Chern character.

This construction works with little change for the more general
twistings~\thetag{2}. 

The differential geometric model we have outlined not only gives geometric
representatives of twisted topological $K$-theory classes, but also geometric
representatives of twisted {\it differential\/} $K$-theory classes.
Similarly, the geometric twistings~\thetag{5} give geometric models for
differential cohomology classes.  See~\cite{HS} for foundations of
(untwisted) differential cohomology theories in general, and~\cite{L} for the
basics differential $K$-theory. 

\head \S3  Loop groups
\endhead

Let $G$~be a compact Lie group.  The loop group~$LG$ of~$G$ is the space of
smooth maps~$S^1\to G$.  There is a twisted description we use instead.
Namely, let $R\to S^1$ be a principal $G$-bundle and $\lg$~the group of gauge
transformations, i.e., the space of smooth sections of the bundle of groups
$G_R\to \cir$ associated to~$R$.  A trivialization of~$R$ gives an
isomorphism~$\lg\cong LG$.  Note that $R$~is necessarily trivializable if
$G$~is connected, and in general the topological class of~$R$ is labeled by a
conjugacy class in~$\pi _0G$.  Let $\GR\subset G$ be the union of components
which comprise that conjugacy class.  The theory of loop groups~$\lg$ is
described in~\cite{PS}, \cite{Ka}, and some specific further developments
appear in~\cite{FF}, \cite{We}.
 
One salient feature of loop groups is the existence of nontrivial central
extensions 
  $$ 1 \longrightarrow \TT \longrightarrow \lgt \longrightarrow
     \lg\longrightarrow 1.  $$
As in~\thetag{2} we may also consider gradings $\epsilon \:\lg\to\zt$.  We
call the pair~$\tau =(\lgt,\epsilon )$ a {\it graded central extension} and
denote it simply as~$\lgt$.  It has an invariant in
  $$ H^1_G\bigl(\GR;\zt\bigr)\times H^3_G\bigl(\GR;\ZZ\bigr) \tag{7} $$
as follows.  Fix a basepoint~$s\in S$ and consider the product
  $$ P=\Cal{A}_R\times R_s \tag{8} $$
of the space of connections~ $\Cal{A}_R$ on~$R$ and the fiber~$R_s$ at the
basepoint.  Then $\lg$~acts freely on~$P$ with quotient the holonomy map
  $$ \hol\:P\longrightarrow \GR.  $$
Furthermore, the $G$~action on~$\GR$ by conjugation lifts to the $G$~action
on the $R_s$~factor of~$P$, where it commutes with the $\lg$~action.  Thus we
have an equivariant twisting 
  $$ \ttil = (\lg,\epsilon ,\lgt,P) $$
whose isomorphism class, called the {\it level\/}, lies in~\thetag{7}.  (The
basepoint is not necessary.  It is more natural to replace~\thetag{8} by
$\Cal{A}_R\times R$ and $\lg$~by the group of bundle automorphisms which
cover rotations in the base~$S^1$.)
 
As  a  warm-up to  representations  of  loop  groups, recall  the  Borel-Weil
construction of  representations of a  compact connected Lie  group~$G$.  Let
$T\subset G$ be a maximal torus, and $F=G/T$~the {\it flag manifold\/}.  Then
$F$~admits $G$-invariant complex  structures.  Fix one.  A character~$\lambda
$  of~$T$ determines  a holomorphic  line bundle~$L_\lambda  \to F$,  and the
standard construction  takes the induced virtual representation  of~$G$ to be
$\oplus _{q}(-1)^qH^q(F;L_\lambda  )$, which is the  $G$-equivariant index of
the  $\dbar$~operator.   We  can  use  the  Dirac  operator  instead  of  the
$\dbar$~operator.  This has  the advantage that no complex  structure need be
chosen, though  it can be  described in holomorphic  terms as the  {\it $\rho
$-shift\/}  $ L_\lambda \longrightarrow  L_\lambda \otimes  K^{1/2}_F$, where
$K^{1/2}_F\to F$ is  a square root of the  canonical bundle.  More generally,
if $Z\subset G$  is the centralizer of any subtorus  of~$T$, and $F'=G/Z$ the
corresponding generalized flag manifold, then there is a Dirac induction map
  $$ K_Z\cong K_G(F')\longrightarrow K_G. \tag{9} $$
A representation of~$Z$ defines a $G$-equivariant vector bundle on~$F'$, and
\thetag{9}~is the equivariant index of the Dirac operator with coefficients
in this bundle.
 
A similar construction works for loop groups.  Fix a conjugacy
class~$C\subset \GR$.  The group $\lg\times G$ acts transitively
on~$\hol\inv (C)$.  Let $Z_C$~ denote the stabilizer at some point; it embeds
isomorphically into either factor of~$\lg\times G$.  Introduce the flag
manifold~$\fc=\hol\inv (C)/G$; then the loop group~$\lg$ acts
transitively~$\fc$ with stabilizer the image of~$Z_C$ in~$\lg$.  On the
other hand, the image of~$Z_C$ in~$G$ is the centralizer of an element
in~$C$.  The geometry of~$\fc$ is similar to that of finite dimensional flag
manifolds~\cite{F4}.  An important special case is~$C=\{1\}$ and $R\to\cir$
trivial.  Then the flag manifold is the ``loop grassmannian''~$\scrF=LG/G$.

There is a Dirac induction for loop groups which generalizes~\thetag{9}.  For
this we introduce spinors on~$\lg$ and a canonical graded central extension.
Consider the principal $SO(\frak{g}\oplus \RR)$-bundle associated to
$R\to\cir$ via the twisted adjoint homomorphism
  $$ \aligned
      G&\longrightarrow SO(\frak{g}\oplus \RR) \\
      g&\longmapsto Ad_g\oplus \det Ad_g\endaligned  $$
It is trivializable, and any trivialization induces an isomorphism~$\lg\to
LSO(N)$ for $N=\dim (G)+1$.  There is a canonical graded central extension
of~$LSO(N)$ from the spin representation, and it pulls back to the desired
canonical graded central extension ~$\lgc\to \lg$.  The spin representation
itself defines (projective) spinors on~$\lg$.  For a conjugacy class
$C\subset G(C)$ the embedding $i_C\: Z_C\hookrightarrow \lg$ induces a
graded central extension~$Z_C^{i_C^*\sigma _R(G)}$ which may be described
instead by the real adjoint representation~$\frak{z}_C$, viewed as a
$Z_C$-equivariant vector bundle over a point (see~\thetag{3}).  In this form
it carries a degree shift, and so it is natural to also include a degree
shift in~$\sigma _R(G)$.  We do not specify it precisely, but remark that its
parity agrees with that of~$\dim Z_C$.  If $G$~is simple, connected, and
simply connected and $R\to\cir$ is trivial, then ~$\sigma _R(G)$ has degree
shift~$\dim G$ and the $H^3$~component of the level in~\thetag{7} is the dual
Coxeter number of~$G$ times a generator; the $H^1$~component vanishes.

For a fixed graded central extension~$\lgt$ of~$\lg$ there is a finite set
of isomorphism classes of irreducible positive energy representations
of~$\lgt$ on which the center acts by scalar multiplication.  Let $R^\tau
(G)$~denote the abelian group generated by these equivalence classes.  It is
natural to extend this to a $\ZZ$-graded group with mod~2 periodicity and
possibly nontrivial groups in odd degrees.  Now we describe Dirac induction
for loop groups.  As in the first map of~\thetag{9} a representation of~$Z_C$
defines an $\lg$-equivariant vector bundle over the flag manifold~$\fc$.
However, we are interested in $\lgt$-equivariant vector bundles, so need to
start with a representation of the central extension of~$Z_C$ defined
by~$i_C^*\tau $.  Finally, spinors on the flag manifold~$\fc$ may be
constructed from spinors on~$\lg$ by ``subtracting'' spinors on the adjoint
representation of~$Z_C$, and this imposes an additional twisting.
Altogether, then, Dirac induction is a map
  $$ K_{\zc}^{i_C^*\tau -\frak{z}_C} \longrightarrow R^{\tau -\sigma
     _R(G)}\bigl(\GR\bigr). \tag{10} $$
The {\it adjoint shift\/}~$\sigma _R(G)$ on the right hand side means that we
obtain representations of the fiber product of~$\lgt\to\lg$ with the inverse
of~$\lgc\to\lg$, including a degree shift.  For connected, simply
connected~$G$ it suffices to consider $C=\{1\}$, since then \thetag{10}~ is
surjective, but this is not true in general.
 
The inclusion $\tilde i_C\:C\hookrightarrow \GR$ induces a pushforward in
twisted $K$-theory (cf.~\thetag{4}):
  $$ K_{\zc}^{i_C^*\tau -\frak{z}_C}\cong K_G^{\tilde i_C^*\ttil
     +TC-TG|_ C}(C)\longrightarrow K^{\ttil}_G\bigl(\GR\bigr). \tag{11} $$
The maps~\thetag{10} and~\thetag{11} give, for each conjugacy class~$C$, a
correspondence between certain representations of the loop group and a
twisted $K$-theory group.  Our main result is

        \proclaim{Theorem 1}
 These correspondences induce an isomorphism of abelian groups
  $$ R^{\tau -\sigma_R(G)}\bigl(G(R)\bigr) \longrightarrow K^{\ttil
     }_G\bigl(G(R)\bigr). \tag{12} $$
        \endproclaim

There is a transgression from~$H^4(BG;\ZZ)$ to levels~\thetag{7} with
trivial first component.  More generally, there is an extension
  $$ 0 \longrightarrow H^4(BG;\ZZ) \longrightarrow E^4(BG) \longrightarrow
     H^2(BG;\zt) \longrightarrow 0,  $$
and elements of~$E^4(BG)$ transgress to general levels for all loop groups~
$\lg$ simultaneously.  If $G$~is connected then any $R\to\cir$ is
trivializable, and both sides of~\thetag{12} have a ring structure for any
fixed~$R$.  For arbitrary~$G$ the sum of each side of~\thetag{12} over
representatives of each topological type of $R\to\cir$ has a ring structure.
The multiplication on representations is the {\it fusion product\/}~\cite{V},
\cite{F}, \cite{T}; on twisted $K$-theory it is the pushforward by
multiplication $G\times G\to G$ or equivalently the Pontrjagin product in
$K$-homology.

        \proclaim{Theorem 2}
 If the level of the twisting $\ttil + \tau (TG)$ is transgressed, then the
isomorphisms~ \thetag{12} are compatible with the ring structure.
        \endproclaim

Our proof reduces both sides of~\thetag{12} to the statement for tori, where
there is a direct argument.

\widestnumber\key{SSSSSSSSSS}   

\enddocument